\catcode`\@=11

\font\tenmsa=msam10
\font\sevenmsa=msam7
\font\fivemsa=msam5
\font\tenmsb=msbm10
\font\sevenmsb=msbm7
\font\fivemsb=msbm5
\newfam\msafam
\newfam\msbfam
\textfont\msafam=\tenmsa  \scriptfont\msafam=\sevenmsa
  \scriptscriptfont\msafam=\fivemsa
\textfont\msbfam=\tenmsb  \scriptfont\msbfam=\sevenmsb
  \scriptscriptfont\msbfam=\fivemsb

\def\hexnumber@#1{\ifnum#1<10 \number#1\else
 \ifnum#1=10 A\else\ifnum#1=11 B\else\ifnum#1=12 C\else
 \ifnum#1=13 D\else\ifnum#1=14 E\else\ifnum#1=15 F\fi\fi\fi\fi\fi\fi\fi}

\def\msa@{\hexnumber@\msafam}
\def\msb@{\hexnumber@\msbfam}
\mathchardef\boxdot="2\msa@00
\mathchardef\boxplus="2\msa@01
\mathchardef\boxtimes="2\msa@02
\mathchardef\square="0\msa@03
\mathchardef\blacksquare="0\msa@04
\mathchardef\centerdot="2\msa@05
\mathchardef\lozenge="0\msa@06
\mathchardef\blacklozenge="0\msa@07
\mathchardef\circlearrowright="3\msa@08
\mathchardef\circlearrowleft="3\msa@09
\mathchardef\rightleftharpoons="3\msa@0A
\mathchardef\leftrightharpoons="3\msa@0B
\mathchardef\boxminus="2\msa@0C
\mathchardef\Vdash="3\msa@0D
\mathchardef\Vvdash="3\msa@0E
\mathchardef\vDash="3\msa@0F
\mathchardef\twoheadrightarrow="3\msa@10
\mathchardef\twoheadleftarrow="3\msa@11
\mathchardef\leftleftarrows="3\msa@12
\mathchardef\rightrightarrows="3\msa@13
\mathchardef\upuparrows="3\msa@14
\mathchardef\downdownarrows="3\msa@15
\mathchardef\upharpoonright="3\msa@16

\mathchardef\downharpoonright="3\msa@17
\mathchardef\upharpoonleft="3\msa@18
\mathchardef\downharpoonleft="3\msa@19
\mathchardef\rightarrowtail="3\msa@1A
\mathchardef\leftarrowtail="3\msa@1B
\mathchardef\leftrightarrows="3\msa@1C
\mathchardef\rightleftarrows="3\msa@1D
\mathchardef\Lsh="3\msa@1E
\mathchardef\Rsh="3\msa@1F
\mathchardef\rightsquigarrow="3\msa@20
\mathchardef\leftrightsquigarrow="3\msa@21
\mathchardef\looparrowleft="3\msa@22
\mathchardef\looparrowright="3\msa@23
\mathchardef\circeq="3\msa@24
\mathchardef\succsim="3\msa@25
\mathchardef\gtrsim="3\msa@26
\mathchardef\gtrapprox="3\msa@27
\mathchardef\multimap="3\msa@28
\mathchardef\therefore="3\msa@29
\mathchardef\because="3\msa@2A
\mathchardef\doteqdot="3\msa@2B

\mathchardef\triangleq="3\msa@2C
\mathchardef\precsim="3\msa@2D
\mathchardef\lesssim="3\msa@2E
\mathchardef\lessapprox="3\msa@2F
\mathchardef\eqslantless="3\msa@30
\mathchardef\eqslantgtr="3\msa@31
\mathchardef\curlyeqprec="3\msa@32
\mathchardef\curlyeqsucc="3\msa@33
\mathchardef\preccurlyeq="3\msa@34
\mathchardef\leqq="3\msa@35
\mathchardef\leqslant="3\msa@36
\mathchardef\lessgtr="3\msa@37
\mathchardef\backprime="0\msa@38
\mathchardef\risingdotseq="3\msa@3A
\mathchardef\fallingdotseq="3\msa@3B
\mathchardef\succcurlyeq="3\msa@3C
\mathchardef\geqq="3\msa@3D
\mathchardef\geqslant="3\msa@3E
\mathchardef\gtrless="3\msa@3F
\mathchardef\sqsubset="3\msa@40
\mathchardef\sqsupset="3\msa@41
\mathchardef\trianglerighteq="3\msa@44
\mathchardef\trianglelefteq="3\msa@45
\mathchardef\bigstar="0\msa@46
\mathchardef\between="3\msa@47
\mathchardef\blacktriangledown="0\msa@48
\mathchardef\blacktriangleright="3\msa@49
\mathchardef\blacktriangleleft="3\msa@4A
\mathchardef\blacktriangle="0\msa@4E
\mathchardef\triangledown="0\msa@4F
\mathchardef\eqcirc="3\msa@50
\mathchardef\lesseqgtr="3\msa@51
\mathchardef\gtreqless="3\msa@52
\mathchardef\lesseqqgtr="3\msa@53
\mathchardef\gtreqqless="3\msa@54
\mathchardef\Rrightarrow="3\msa@56
\mathchardef\Lleftarrow="3\msa@57
\mathchardef\veebar="2\msa@59
\mathchardef\barwedge="2\msa@5A
\mathchardef\doublebarwedge="2\msa@5B
\mathchardef\angle="0\msa@5C
\mathchardef\measuredangle="0\msa@5D
\mathchardef\sphericalangle="0\msa@5E
\mathchardef\varpropto="3\msa@5F
\mathchardef\smallsmile="3\msa@60
\mathchardef\smallfrown="3\msa@61
\mathchardef\Subset="3\msa@62
\mathchardef\Supset="3\msa@63
\mathchardef\Cup="2\msa@64

\mathchardef\Cap="2\msa@65

\mathchardef\curlywedge="2\msa@66
\mathchardef\curlyvee="2\msa@67
\mathchardef\leftthreetimes="2\msa@68
\mathchardef\rightthreetimes="2\msa@69
\mathchardef\subseteqq="3\msa@6A
\mathchardef\supseteqq="3\msa@6B
\mathchardef\bumpeq="3\msa@6C
\mathchardef\Bumpeq="3\msa@6D
\mathchardef\lll="3\msa@6E

\mathchardef\ggg="3\msa@6F

\mathchardef\circledS="0\msa@73
\mathchardef\pitchfork="3\msa@74
\mathchardef\dotplus="2\msa@75
\mathchardef\backsim="3\msa@76
\mathchardef\backsimeq="3\msa@77
\mathchardef\complement="0\msa@7B
\mathchardef\intercal="2\msa@7C
\mathchardef\circledcirc="2\msa@7D
\mathchardef\circledast="2\msa@7E
\mathchardef\circleddash="2\msa@7F
\def\ulcorner{\delimiter"4\msa@70\msa@70 }
\def\urcorner{\delimiter"5\msa@71\msa@71 }
\def\llcorner{\delimiter"4\msa@78\msa@78 }
\def\lrcorner{\delimiter"5\msa@79\msa@79 }
\def\yen{\mathhexbox\msa@55 }
\def\checkmark{\mathhexbox\msa@58 }
\def\circledR{\mathhexbox\msa@72 }
\def\maltese{\mathhexbox\msa@7A }
\mathchardef\lvertneqq="3\msb@00
\mathchardef\gvertneqq="3\msb@01
\mathchardef\nleq="3\msb@02
\mathchardef\ngeq="3\msb@03
\mathchardef\nless="3\msb@04
\mathchardef\ngtr="3\msb@05
\mathchardef\nprec="3\msb@06
\mathchardef\nsucc="3\msb@07
\mathchardef\lneqq="3\msb@08
\mathchardef\gneqq="3\msb@09
\mathchardef\nleqslant="3\msb@0A
\mathchardef\ngeqslant="3\msb@0B
\mathchardef\lneq="3\msb@0C
\mathchardef\gneq="3\msb@0D
\mathchardef\npreceq="3\msb@0E
\mathchardef\nsucceq="3\msb@0F
\mathchardef\precnsim="3\msb@10
\mathchardef\succnsim="3\msb@11
\mathchardef\lnsim="3\msb@12
\mathchardef\gnsim="3\msb@13
\mathchardef\nleqq="3\msb@14
\mathchardef\ngeqq="3\msb@15
\mathchardef\precneqq="3\msb@16
\mathchardef\succneqq="3\msb@17
\mathchardef\precnapprox="3\msb@18
\mathchardef\succnapprox="3\msb@19
\mathchardef\lnapprox="3\msb@1A
\mathchardef\gnapprox="3\msb@1B
\mathchardef\nsim="3\msb@1C
\mathchardef\napprox="3\msb@1D
\mathchardef\nsubseteqq="3\msb@22
\mathchardef\nsupseteqq="3\msb@23
\mathchardef\subsetneqq="3\msb@24
\mathchardef\supsetneqq="3\msb@25
\mathchardef\subsetneq="3\msb@28
\mathchardef\supsetneq="3\msb@29
\mathchardef\nsubseteq="3\msb@2A
\mathchardef\nsupseteq="3\msb@2B
\mathchardef\nparallel="3\msb@2C
\mathchardef\nmid="3\msb@2D
\mathchardef\nshortmid="3\msb@2E
\mathchardef\nshortparallel="3\msb@2F
\mathchardef\nvdash="3\msb@30
\mathchardef\nVdash="3\msb@31
\mathchardef\nvDash="3\msb@32
\mathchardef\nVDash="3\msb@33
\mathchardef\ntrianglerighteq="3\msb@34
\mathchardef\ntrianglelefteq="3\msb@35
\mathchardef\ntriangleleft="3\msb@36
\mathchardef\ntriangleright="3\msb@37
\mathchardef\nleftarrow="3\msb@38
\mathchardef\nrightarrow="3\msb@39
\mathchardef\nLeftarrow="3\msb@3A
\mathchardef\nRightarrow="3\msb@3B
\mathchardef\nLeftrightarrow="3\msb@3C
\mathchardef\nleftrightarrow="3\msb@3D
\mathchardef\divideontimes="2\msb@3E
\mathchardef\varnothing="0\msb@3F
\mathchardef\nexists="0\msb@40
\mathchardef\mho="0\msb@66
\mathchardef\thorn="0\msb@67
\mathchardef\beth="0\msb@69
\mathchardef\gimel="0\msb@6A
\mathchardef\daleth="0\msb@6B
\mathchardef\lessdot="3\msb@6C
\mathchardef\gtrdot="3\msb@6D
\mathchardef\ltimes="2\msb@6E
\mathchardef\rtimes="2\msb@6F
\mathchardef\shortmid="3\msb@70
\mathchardef\shortparallel="3\msb@71
\mathchardef\smallsetminus="2\msb@72
\mathchardef\thicksim="3\msb@73
\mathchardef\thickapprox="3\msb@74
\mathchardef\approxeq="3\msb@75
\mathchardef\succapprox="3\msb@76
\mathchardef\precapprox="3\msb@77
\mathchardef\curvearrowleft="3\msb@78
\mathchardef\curvearrowright="3\msb@79
\mathchardef\digamma="0\msb@7A
\mathchardef\varkappa="0\msb@7B
\mathchardef\hslash="0\msb@7D
\mathchardef\hbar="0\msb@7E
\mathchardef\backepsilon="3\msb@7F
\def\Bbb{\ifmmode\let\next\Bbb@\else
 \def\next{\errmessage{Use \string\Bbb\space only in math mode}}\fi\next}
\def\Bbb@#1{{\Bbb@@{#1}}}
\def\Bbb@@#1{\fam\msbfam#1}

\catcode`\@=\active

\def\inv{^{\raise.15ex\hbox{${
  \scriptscriptstyle -}$}\kern-.05em 1}}

\def\Dsl{\,\raise.15ex\hbox{$/$}\mkern-13.5mu D}
\def\dsl{\raise.15ex\hbox{$/$}\kern-.57em\hbox{$\partial$}}

\def\lspace{\ifx\answ\bigans{}\else\qquad\fi}
\def\del{\partial}

\def\lform{\hbox{$\sqcup$}\llap{\hbox{$\sqcap$}}}
\def\darr#1{\raise1.5ex\hbox{$\leftrightarrow$}
\mkern-16.5mu #1}

\def\INT{{\textstyle \int\kern-.642em\int}}

\def\Z{{\Bbb Z}}

\def\rbiprod{{\cdot\kern-.33em\triangleright\!\!\!<}}
\def\lbiprod{{>\!\!\!\triangleleft\kern-.33em\cdot}}

\def\tens{\mathop{\otimes}}

\def\<{\langle}
\def\>{\rangle}

\def\vecx{{\bf x}}
\def\vecl{{\bf l}}

\def\<{\langle}
\def\>{\rangle}

\def\equad{\kern -1.7em}
\def\qqquad{\qquad\quad}
\def\nquad{{\!\!\!\!\!\!}}

\def\und#1{{\underline {#1}}}

\def\text#1{\mbox{\rm #1}}
\def\note#1{}

\def\blacksquare{{\lform}}
\def\frac#1#2{{{#1\over#2}}}

\def\proof{\goodbreak\noindent{\bf Proof\quad}}

\def\endproof{{\ $\lform$}\bigskip }

\def\eqn#1#2{\begin{equation}#2\label{#1}\end{equation}}

\def\align#1{\begin{eqnarray*}#1\end{eqnarray*}}

\def\ceqn#1#2{\begin{equation}\label{#1}\begin{array}{c}#2\end{array}
\end{equation}}

\documentstyle[11pt]{article}
\newtheorem{lemma}{Lemma}[section] \newtheorem{propos}[lemma]{Proposition}
 \newtheorem{theorem}[lemma]{Theorem}

\begin{document}
\baselineskip 22pt
{\ }\qquad\qquad \hskip 4.3in DAMTP/95-58
\vspace{.2in}

\begin{center} {\LARGE  FERMIONIC $q$-FOCK SPACE AND BRAIDED GEOMETRY}
\\ \baselineskip 13pt{\ }
{\ }\\ S. Majid\footnote{Royal Society University Research Fellow and Fellow of
Pembroke College, Cambridge. On leave  1995 + 1996 at the Department of
Mathematics, Harvard University, Cambridge MA02138, USA}\\
{\ }\\
Department of Applied Mathematics \& Theoretical Physics\\
University of Cambridge, Cambridge CB3 9EW\\
+\\
Research Institute of Mathematical Sciences\\
Kyoto University, Kyoto 606, Japan
\end{center}
\begin{center}
July -- revised November, 1995
\end{center}

\vspace{10pt}
\begin{quote}\baselineskip 13pt
\noindent{\bf Abstract}
We write the fermionic $q$-Fock space representation of $U_q(\hat{sl_n})$ as
an infinite extended braided tensor product of finite-dimensional fermionic
$U_q(sl_n)$-quantum planes or exterior algebras. Using braided geometrical
techniques developed for such
quantum exterior algebras,
we provide a new approach to the Kashiwara-Miwa-Stern  action of the Heisenberg
algebra on the $q$-fermionic Fock space, obtaining the action in detail for the
lowest nontrivial case $[b_{2},b_{-2}]=2({1-q^{-4n}\over 1-q^{-4}})$.
Our R-matrix approach includes other Hecke R-matrices as well.


\bigskip

\end{quote}
\baselineskip 22pt

\section{Introduction}

In this note we use techniques from `braided geometry' to study the
$q$-deformed fermionic Fock space representations of
the affine quantum groups $U_q(\hat{sl_n})$
\cite{Hya}\cite{MisMiw}\cite{Ste}\cite{KasMiwSte}. The properties of this
$q$-deformed Fock space are closely connected with the theory of vertex
operator algebras and $q$-correlation
functions. In particular, using the vertex operator algebra approach it has
been shown in \cite{KasMiwSte} that there is an
action of the Heisenberg algebra on the level 1 fermionic Fock space
representation of $U_q(\hat{sl_n})$ through natural `shift' operators $b_i$.

We provide now a new approach to this $q$-fermionic Fock space via the theory
of  braided groups\cite{Ma:introp} as developed extensively by the author in
recent years. We refer to
\cite{Ma:varen} for a more recent review.  The
standard finite-dimensional quantum planes have such a braided group
structure or coaddition, which allows
one to define braided differentiation\cite{Ma:fre}, integration, epsilon
tensors\cite{Ma:eps}, differential forms, etc. on such spaces in a systematic
way. Using
such techniques, we explicitly derive the Heisenberg algebra action of
\cite{KasMiwSte} for the lowest
non-trivial  generators $b_1,b_2$. Even these cases will be hard enough, but we
believe that they
demonstrate the possibility of a new approach using such techniques. Ultimately
it may be possible to compute
$q$-correlation functions themselves by such methods, which is one of the
motivations for the work.

Our starting point is the infinite-dimensional quantum
planes or exchange algebras, associated
to unitary solutions of the parametriced Yang-Baxter equations
\eqn{pybe}{ R_{12}({z\over
w})R_{13}(z)R_{23}(w)=R_{23}(w)R_{13}(z)R_{12}({z\over w}),\quad
R(z)=R(z^{-1})_{21}^{-1}}
in a compact notation. Associated to this is the corresponding fermionic
quantum plane $\Lambda(R(z))$ with
\eqn{fock}{\theta_1(z)\theta_2(w)=-\theta_2(w)\theta_1(z) R({z\over w}),\quad
{\rm i.e.}\quad
\theta_i(z)\theta_j(w)=\theta_b(w)\theta_a(z) R^a{}_i{}^b{}_j({z\over w})}
where $R(z)\in M_n\tens M_n$ and $\theta(z)_i$, $i=1,\cdots ,n$. There are also
similar formulae without
the - signs, for bosonic-type exchange algebras.  The fermionic Fock
space
in \cite{KasMiwSte} is of this general type (\ref{fock}), where, more
precisely, the authors considered vector
near to a chosen `vacuum vector', rather than the algebra itself. We refer to
\cite{KasMiwSte} for details on this final step.

In Section  2, we study the algebra (\ref{fock}) for the entire class of
solutions of (\ref{pybe}) of the
form
\eqn{baxt}{ R(z)={R-zR_{21}^{-1}\over q-zq^{-1}}.}
This Baxterisation formula solves (\ref{pybe}) for {\em any} matrix solution
$R$ of the ordinary Yang-Baxter
equations which is of Hecke type, in the sense
\eqn{hecke}{ (PR-q)(PR+q^{-1})=0,}
where $P$ is the permutation matrix, which is the generality at which we work.
This approach
 includes the $U_q(\hat{sl_n})$ R-matrix
as well as other more nonstandard systems. We show that the algebra
$\Lambda(R(z))$ is an infinite `tensor product' of copies   of the fermionic
quantum plane $\Lambda(R)$ with
\eqn{ext}{ \theta_1\theta_2=-q\theta_2\theta_1 R.}
Such fermionic quantum planes have key properties from the theory of braided
geometry, which we shall use.
Among them is the braided coaddition
\eqn{coadd}{ \Delta\theta=\theta\tens 1+1\tens\theta,\quad
(1\tens\theta_1)(\theta_2\tens 1)
=-q^{-1}(\theta_2\tens 1)(1\tens\theta_1)R}
where the two copies of $\Lambda(R)$ in $\Lambda(R)\und\tens \Lambda(R)$ enjoy
the braid
statistics shown (generalising the usual Bose-Fermi statistics of usual
exterior algebras), which
makes them braided groups rather than quantum groups. Moreover,
 because braided geometry works as well for fermionic
as for bosonic spaces, its principal notions such as braided-differentiation,
etc.,  work as well for  $\Lambda(R)$ as for the
more usual bosonic quantum planes. In particular, as a case of \cite{Ma:fre},
we have braided differentiation
on fermionic quantum planes\cite{Ma:eps}
\ceqn{bradif}{
\del^i(\theta_1\theta_2\cdots\theta_m)=e_1^i\theta_2\cdots\theta_m
[m,-q^{-1}R]_{1\cdots m}\\
 \theta_1\theta_2\cdots\theta_m\overleftarrow{\del^i}
=\theta_1\cdots\theta_{m-1}e^i_m
\overline{[m; -q^{-1}R]}_{1\cdots m}\\
{}[m,R]_{1\cdots m}=1+(PR)_{12}+(PR)_{12}(PR)_{23}+\cdots +(PR)_{12}\cdots
(PR)_{m-1m}\\
{}\overline{[m;R]}_{1\cdots m}=1+(PR)_{m-1m}+(PR)_{m-1m}(PR)_{m-2m-1}+\cdots +
(PR)_{m-1m}\cdots (PR)_{12}}
as operators $\del^i,\overleftarrow{\del^i}:\Lambda(R)\to \Lambda(R)$. Here
$(e^i)_j=\delta^i{}_j$ is a basis vector.
One can also apply such ideas at the
infinite-dimensional level (\ref{fock}), as functional
differentiation, though we do not do so here.

Our goal is to make use of some of the rich structure of finite-dimensional
braided spaces to study the infinite-dimensional
fermionic Fock space. In effect, we
study these exchange algebras as `braided wave functions' where at each point
(in momentum space) we have a mode $\theta^{i}$
behaving as a fermionic quantum plane. Moreover, our deriviations in this paper
do not depend at any point on the
precise form of the Hecke R-matrix. Hence we include not only the
$U_q(\hat{sl_n})$
theory but, in principle, generalise it to other non-standard affine quanutm
groups associated to the Baxterisation (\ref{baxt})
of other Hecke R-matrices as well. We derive the Heisenberg algebra action in
Section~3
in this setting. In Section~4, we conclude with some comments
about covariance.

Some notations in the paper are as follows. Apart from the {\em braided integer
matrices}\cite{Ma:fre} $[m,R]$ and $\overline{[m,R]}$
in (\ref{bradif}), we also set
\[[m;q^{-2}]\equiv{1-q^{-2m}\over 1-q^{-2}},\quad
[m,n;R]\equiv(PR)_{mm+1}(PR)_{m+1m+2}\cdots (PR)_{n-1n}\]
\[\overline{[m,n;R]}\equiv(PR)_{n-1n}\cdots (PR)_{m+1m+2}(PR)_{mm+1}.\]
There is a change in conventions $q\to q^{-1}$ in our paper relative to
\cite{KasMiwSte}. Also, we write the fermionic quantum
plane relations such as (\ref{ext}) in the even more compact form in which we
suppress the numerical suffices entirely. Thus
\[ \theta\theta\equiv\theta_1\theta_2,\quad {\rm i.e.},\quad
\theta\theta=-q\theta\theta PR\]
is (\ref{ext}) in our notation: the tensor product of the vector indices
$\theta_i$ is to be understood. When we do write
numerical suffices $\theta_1,\theta_2$ etc, we henceforth mean the actual
components of the vector $\theta$. Finally, we
write
\[ \{\theta,\psi\}_R\equiv \theta\psi+q^{-1}\psi\theta PR,\quad{\rm i.e.}\quad
\{\theta_i,\psi_j\}_R
\equiv \theta_i\psi_j+q^{-1}\psi_b\theta_a R^a{}_i{}^b{}_j\]
and sometimes ${\bf R}\equiv -q^{-1}R$, as  useful shorthand notations.

\subsection*{Acknowledgements} These results were obtained during a visit in
June 1995 to R.I.M.S. in Kyoto under a joint programme with the Isaac Newton
Institute in Cambridge and the J.S.P.S. I would like to thank my host T. Miwa
for extensive discussions.

\section{Fermionic Fock space}

The level 1 Fock space representation of $U_q(\hat{sl_n})$ has been constructed
in \cite{Hya}\cite{MisMiw} and studied
further in several papers, notably \cite{Ste}\cite{KasMiwSte}. Here we take a
slightly
different point of view on this
representation, taking as starting point the fermionic `exchange algebra'
$\Lambda(R(z))$ defined in (\ref{fock}).
 Our goal in this section is to break down the structure of
this exchange algebra into many copies of standard finite-dimensional fermionic
quantum planes $\Lambda(R)$ as in
(\ref{ext}). We write $\theta(z)=\sum_{z\in \Z}\theta^{i}z^i$.

\begin{theorem} When $R(z)$ is of the form (\ref{baxt}) (as in the $sl_n$ case)
then $\Lambda(R(z))$ is an
infinite number of copies $\{\theta^{i}\}$ of the fermionic quantum plane
$\Lambda(R)$ associated to the
finite-dimensional R-matrix $R$, with relations
\[ \theta^{i}\theta^{i}(PR+q^{-1})=0,\quad \{\theta^{i},\theta^{i-1}\}_R=0\]
\[ \{\theta^{i},\theta^{j}\}_R=(q^{-2}-1)\left(\sum_{s=1}^{s<{i-j\over
2}}\theta^{j+s}
\theta^{i-s}(1+q^{-2})^{s-1}(1+P{\bf R})+\theta^{{i+j\over 2}}\theta^{{i+j\over
2}}q^{-2({i-j\over 2}-1)}\right)\]
for $i-j>1$. Here the last term is included only if $i-j$ is even.
\end{theorem}
\proof From the form of $R(z)$ we have
\[ \sum_{i,j}(q-q^{-1}{z\over w})\theta^{i}\theta^{j} z^i w^j
=\sum_{i,j}\theta^{j}w^j\theta^{i}z^i (PR-{z\over w}(PR)^{-1}).\]
We equate powers of $z,w$, and hence
require
\eqn{moderel}{ \theta^{j}\theta^{i} PR + \theta^{i}\theta^{j} q=
\theta^{j+1}\theta^{i-1} (PR)^{-1}+\theta^{i-1}
\theta^{j+1} q^{-1}.}
Considering the same equation with $i\to j+1$ and $j\to i-1$ and combining with
(\ref{moderel}) times $qPR$, gives
\eqn{modeanticom}{ (\theta^{i}\theta^{j}+\theta^{j}\theta^{i})(PR+q^{-1})=0,}
on using the Hecke condition (\ref{hecke}). This implies, in particular, that
the $\theta^{i}$ modes each obey the finite-dimensional fermionic quantum plane
algebra. Next, we consider (\ref{moderel}) with $j=i-1$, i.e.,
\[ \theta^{i-1}\theta^{i}PR+\theta^{i}\theta^{i-1} q=\theta^{i}\theta^{i-1}
(PR)^{-1}+
\theta^{i-1}\theta^{i} q^{-1}.\]
Combining with (\ref{modeanticom}) and the Hecke condition
$(PR)^2=1+(q-q^{-1})PR$  gives $\{\theta^{i},\theta^{i-1}\}_R=0$
for neighbouring modes.
Finally, for non-neighbouring modes, we use the Hecke condition to write
(\ref{moderel}) in the form
\eqn{anticomind}{
\{\theta^{i},\theta^{j}\}_R=\{\theta^{i-1},\theta^{j+1}\}_R+(q^{-2}-1)
(\theta^{i-1}\theta^{j+1}+\theta^{j+1}\theta^{i-1}),}
which gives an inductive formula for $\{\theta^{i},\theta^{j}\}_R$ in terms of
`usual' anticommutators of the intermediate modes. Alternatively, which we
prefer, we use (\ref{modeanticom}) and the
Hecke condition to write (\ref{anticomind}) as
\eqn{modeind}{\{\theta^{i},\theta^{j}\}_R=\{\theta^{i-1},
\theta^{j+1}\}_R(1+(q-q^{-1})PR)+
(q^{-2}-1)\theta^{j+1}\theta^{i-1}(1-q^{-1}PR).}
Using this, we obtain the formula stated for the ordering relations between
non-adjacent modes, by induction. Note that, by the Hecke
condition (\ref{hecke}), $(1-q^{-1}PR)PR=(1-q^{-1}PR)(-q^{-1})$. The start of
the induction
is when the $i,j$ are equal or one apart (as $i-j$ is even or odd), which cases
we have already computed separately.
We see that between adjacent modes there are the usual braid statistics
associated to two
copies of the finite-dimensional fermionic quantum plane (as needed for their
braided coaddition structure in (\ref{coadd})).
Between modes that are further apart, we have the same `leading' braid
statistics + decendent terms involving
intermediate modes. \endproof

The algebra in this theorem is computed formally from the powerseries, but can
afterwards be taken as a definition
of
the exchange algebra, as generated by $\theta^{i}$. We proceed now on this
basis. We see that
each of the modes has a geometrical picture as the algebra $\Lambda(R)$ of
$q$-differential forms;
see \cite{Ma:eps} for the braided-geometrical construction (starting from the
braided coaddition law). In particular,
in nice cases (such as the $sl_n$ case), each has a top form
\[ \omega^{i}=\theta^{i}_1\cdots\theta^{i}_n\]
with all others of this degree being multiplies of it. The products
$\theta^{i}\omega^{i}$ are zero for all $i$.
There is also an underlying bosonic space with $\theta^{i}=d\vecx^{i}$, where
$\vecx^{i}$ obey $\vecx^{i}\vecx^{i}
=\vecx^{i}\vecx^{i}q^{-1}PR$. We do not use this full geometrical picture here,
regarding the $\theta^{i}$ as intrinsic
fermionic-type coordinates in their own right.

It is worth noting that our fermionic Fock space algebra in Theorem~2.1 is
clearly a more complicated variant
of the  actual braided tensor product algebra $\und\tens_{i=-\infty}^{\infty}
\Lambda^{i}(R)$ with relations
\eqn{gerv}{ \theta^{i}\theta^{i}(PR+q^{-1})=0,\quad
\{\theta^{i},\theta^{j}\}_R=0}
for all $i>j$. This algebra was discussed in \cite{Ma:introp}, where it was
proposed as a discrete model of the exhange
algebra in 2-D quantum gravity\cite{Ger}. Indeed, one can consider it as a
fermionic  exchange algebra for the discretely (and additively)
parametrised R-matrix
\eqn{gervR}{ R(i-j)=\cases{q^{-1}R& $i>j$\cr qR&$i=j$\cr qR_{21}^{-1}& $i<j$.}}
The algebra (\ref{gerv}), although pertaining to a different model than the one
above (and with $i$ as a discrete version of
a position variable rather than a mode label), nevertheless has  a similar form
to our fermionic Fock space in Theorem~2.1,
just without the descendent modes.
Moreover, its construction as a braided tensor product (with relations as in
(\ref{coadd}) between different modes)
ensures that it remains covariant under (a dilatonic
extension of) $U_q(sl_n)$ or other quantum group (according to the R-matrix).
By contrast, the more complicated
fermionic Fock space in Theorem~2.1 is covariant under $U_q(\hat{sl_n})$ or
other affine quantum group.

\section{Computation of the Heisenberg algebra action}

It is clear from the form of the relations (\ref{moderel}) of $\Lambda(R(z))$
that
\eqn{bi}{b_i:\Lambda(R(z))\to \Lambda(R(z)),\quad b_i(\theta^{j})=\theta^{j+i}}
is a derivation on the algebra, for each $i$.  It is shown in \cite{KasMiwSte},
(by Hecke algebra and vertex operator methods) that these $b_i$ define an
action of the Heisenberg algebra
according to
\eqn{heis}{ [b_i,b_{-j}]=\delta_{i,j} i\left({1-q^{-2ni}\over
1-q^{-2i}}\right),}
when acting on
\[ \omega=\omega^{0}\omega^{1}\cdots\]
or vectors near to this (differing only in finitely many coefficients).  We
show now how this result can alternatively
be obtained by braided-geometrical methods. Note that $\omega$ is in a
completion of the algebra generated by the modes. However,
all our operations stay within the space of vectors near to it, and hence
remain algebraic; see \cite{KasMiwSte} for a more
formal way to say this.

\begin{propos}  For $i\ge 1$, we have
\[ b_i(\omega)=0,\quad b_{-i}(\omega)=b_{-i}(\omega^0)\omega^1\cdots+\omega^0
b_{-i}(\omega^1)\omega^2\cdots+\cdots + \omega^0\omega^1\cdots
\omega^{i-2} b_{-i}(\omega^{i-1})\omega^i\cdots.\]
\end{propos}
\proof Firstly, $b_i(\omega)=0$ for $i\ge 1$ since $b_i(\omega^j)$ has in it
modes $\theta^{j+i}$; moving these to the right
using the braided-anticommutation relations with
$\theta^j,\theta^{j+1},\cdots,\theta^{j+i-1}$, gives eventually
$\theta^{j+i}\omega^{j+i}=0$.
Along the way, if $i\ge 2$, we generate decendents  which lie in the range
$\theta^{j+1},
\cdots,\theta^{j+i-1}$; moving each of these to the right kills these as well.
Similarly for their descendents, etc.

For $b_{-i}$ we have
\[ b_{-i}(\omega^j)=\theta^{j-i}_1\theta^j_2\cdots\theta_n^j+\cdots
+\theta_1^j\cdots\theta^j_{n-1}\theta_n^{j-i}=\theta^{j-i}_{a_1}
\theta^j_{a_2}\cdots\theta^j_{a_n}[n;{\bf R}]^{a_1\cdots a_n}_{1\cdots n}+{\rm
decsendents}\]
where the decendents involve
 $\theta^{j-i+1},\cdots,\theta^{j-1}$. We moved $\theta^{j-2}$ to the left in
each term, just as in the definition of
braided differentiation\cite{Ma:fre}, but now picking up descendents from the
right hand side of the anticommutators in Theorem~2.1.

Hence, when we compute $b_{-i}(\omega)$ as a derivation, only the first $i$
terms
contribute, as stated; the  $\omega^0\cdots \omega^{j-1}b_{-i}(\omega^j)$ for
$j\ge i$ do not contribute because the terms
of $b_{-i}(\omega^j)$ each contain a mode in the range
$\theta^{j-i},\cdots,\theta^{j-1}$ which, using the
relations in Theorem~2.1,
can be pushed left until it  multiplies one of
$\omega^{j-i},\cdots,\omega^{j-1}$, and
thereby vanishes. The descendents generated in this process when $i\ge 2$ can
likewise be pushed to the left and anihilated.
Similarly for their descendents, etc. \endproof

The simplest case of (\ref{heis}) follows trivially:

\begin{propos}
$b_{-1}(\omega^j)=\theta^{j-1}_{a_1}\theta^{j}_{a_2}\cdots
\theta^{j}_{a_n}[n;{\bf R}]^{a_1\cdots a_n}_{1\cdots n}$. Hence
$[b_1,b_{-1}]=[n,q^{-2}]$
when acting on $\omega$.
\end{propos}
\proof In this case $\theta^{j-1}$ is adjacent to $\theta^{j}$ so no
descendents are generated when we move it to the
left in each term of $b_{-1}(\omega^j)$. Hence
$b_{-1}(\omega)=\theta^{-1}\theta^0\cdots\theta^0[n;{\bf R}]\omega^1\cdots$.
When we apply $b_1$ to this, only the action on $\theta^{-1}$ contributes:
other modes have degree $\ge 1$ and anihilate when
moved to the right. Hence $b_1(b_{-1}(\omega))=\theta^0\cdots\theta^0[n;{\bf
R}]\omega^1\cdots$. On the other hand, $PR$
acts as $-q^{-1}$ on $\theta\theta$ (the defining relations of each mode
$\Lambda(R)$ in Theorem~2.1). Hence $[n;{\bf R}]$
can be replaced by $[n;q^{-2}]$ when acting on $\Lambda^{(0)}(R)$. \endproof

The same techniques apply for the action of the higher Heisenberg generators.
We do the computation now for $[b_2,b_{-2}]$.

\begin{lemma}
\align{b_{-2}(\omega^j)&=&\theta^{j-2}\theta^j\cdots\theta^j[n;{\bf
R}]_{1\cdots n}\\
&&\quad +(q^{-2}-1)\theta^{j-1}\theta^{j-1}\theta^j\cdots\theta^j(
[n-1;{\bf R}]_{2\cdots n}+[2,3;{\bf R}][1,2;{\bf R}]
[n-2;{\bf R}]_{3\cdots n}\\
&&\qqquad\qquad +\cdots+[2,n-1;{\bf R}][1,n-2;{\bf R}][2;{\bf
R}]_{n-1n}+[2,n;{\bf R}][1,n-1;{\bf R}]).}
Hence
\[ b_2(b_{-2}(\omega^0))\omega^1\cdots=\left([n;q^{-2}]+(1-q^{-2})
\left([n-1;q^{-4}]-q^{-2(n-1)}[n-1;q^{-2}]\right)\right)\omega.\]
\end{lemma}
\proof Clearly,
\align{b_{-2}(\omega^j)&=&\theta^{j-2}_1\theta^j_2\cdots\theta^j_n
+\cdots+\theta^j_1\cdots\theta^j_{n-1}\theta^{j-2}_n\\
&=&\theta^{j-2}\theta^j\cdots\theta^j[n,{\bf R}]_{1\cdots
n}+(q^{-2}-1)\theta_1^{j-1}\theta^{j-1}\theta^j\cdots\theta^j
[n-1;{\bf R}]_{2\cdots n}\\
&&\quad
+(q^{-2}-1)\theta_1^j\theta_2^{j-1}\theta^{j-1}\theta^j\cdots\theta^j[n-2;{\bf
R}]_{3\cdots n}+\cdots
+(q^{-2}-1)\theta_1^j\cdots\theta_{n-2}^j\theta_{n-1}^{j-1}\theta_n^{j-1},}
where we use
\[ \theta^j\theta^{j-2}=\theta^{j-2}\theta^j P{\bf R}+
(q^{-2}-1)\theta^{j-1}\theta^{j-1}\]
from Theorem~2.1. We move each $\theta^{j-2}$ to the left at the price of a
factor $P{\bf R}$ and a $\theta^{j-1}\theta^{j-1}$. We then
add up all the descendents as generated in each position.

{}From this expression, the expression stated for $b_{-2}(\omega^j)$ follows at
once: in each of the descendent terms, we move $\theta^{j-1}\theta^{j-1}$
to the left, accumulating powers of $P{\bf R}$ for each one.

Then $b_2(b_{-2}(\omega^0))\omega^1\cdots $ is computed as follows. When we
apply $b_2$,
only its action on the $\theta^{-2}$ mode or the first $\theta^{-1}$ mode in
$b_{-2}(\omega^0)$ can contribute, since the other cases produce modes which
can
be pushed to the right and anihilated, along with their descendents.
The first of these gives $\theta^0\cdots\theta^0[n;{\bf
R}]\omega^1\cdots=\omega [n;q^{-2}]$
by the relations in $\Lambda^{(0)}(R)$. The second case contains
$\theta^1\theta^{-1}\theta^0\theta^0\cdots\theta^0$ where $\theta^1$ can also
be pushed to the right and anihilated. In the process, however, it contributes
a descendent
\[ \theta^0\theta^0\cdots\theta^0(q^{-2}-1)^2\left(
[n-1;{\bf R}]_{2\cdots n}+[2,3;{\bf R}][1,2;{\bf R}][n-2;{\bf R}]_{3\cdots
n}+\cdots+[2,n;{\bf R}][1,n-1;{\bf R}]
\right)\omega^1\cdots.\]
Finally, using the relations in $\Lambda^{(0)}(R)$, we can replace $P{\bf R}$
by $q^{-2}$, giving
\[ (q^{-2}-1)^2\left([n-1;q^{-2}]+q^{-4}[n-2;q^{-2}]+\cdots
q^{-4(n-2)}[1;q^{-2}]\right)\omega\qqquad\]
\[\qqquad=(1-q^{-2})\left([n-1;q^{-4}]-q^{-2(n-1)}[n-1;q^{-2}]\right)\omega\]
as stated. \endproof

By a strictly analogous computation, we have
\align{ b_2(\omega^j)&=&\theta^j\cdots\theta^j\theta^{j+2}\overline{[n;{\bf
R}]}_{1\cdots n}\\
&&\quad +(q^{-2}-1)\theta^j\cdots\theta^j\theta^{j+1}\theta^{j+1}
(\overline{[n-1;{\bf R}]}_{1\cdots n-1}+\overline{[1,2;{\bf R}]}\,
\overline{[2,3;{\bf R}]}\, \overline{[n-2;{\bf R}]}_{1\cdots n-2}\\
&&\qqquad\qquad +\cdots+\overline{[1,n-2;{\bf R}]}\, \overline{[2,n-1;{\bf
R}]}\, \overline{[2;{\bf R}]}_{12}+
\overline{[1,n-1;{\bf R}]}\, \overline{[2,n;{\bf R}]}),}
showing its descendents explicitly. Here we moved $\theta^{j+2}$ to the right,
and the resulting descendents also to the right.

\begin{propos} $[b_2,b_{-2}]=2\left({1-q^{-4n}\over 1-q^{-4}}\right)$ when
acting on $\omega$.
\end{propos}
\proof We are now ready to compute
\[ b_2(b_{-2}(\omega))=b_2(b_{-2}(\omega^0)\omega^1+\omega^0
b_{-2}(\omega^1))\omega^2\cdots\]
where $b_2(\omega^2)\omega^3$ etc., do not contribute, as in Proposition~3.1
(shifted down by translation invariance).
The first term is the same as $b_2(b_{-2}(\omega^0))\omega^1\cdots$ (for the
same reason) and was computed in Lemma~3.3. The
second term is
\align{b_2(\omega^0b_{-2}(\omega^1))\omega^2\cdots&
=&b_2(\theta^0\cdots\theta^0\theta^{-1}_{a_1}\theta^1_{a_2}\cdots\theta^1_{a_n}
[n;{\bf R}]^{a_1\cdots a_n}_{1\cdots n})\omega^2\cdots\\
&&=b_2(\theta^{-1}\theta^0\cdots\theta^0[1,n+1;{\bf R}]_{1\cdots n a_1}
\theta^1_{a_2}\cdots\theta^1_{a_n}[n;{\bf R}]^{a_1\cdots a_n}_{1\cdots
n})\omega^2\cdots\\
&&=\theta^1\theta^0\cdots\theta^0[1,n+1;{\bf R}]_{1\cdots n a_1}
\theta^1_{a_2}\cdots\theta^1_{a_n}[n;{\bf R}]^{a_1\cdots a_n}_{1\cdots
n}\omega^2\cdots\\
&&=\theta^0\cdots\theta^0\theta^1\overline{[1,n+1;{\bf R}]}[1,n+1;{\bf
R}]_{1\cdots n a_1}
\theta^1_{a_2}\cdots\theta^1_{a_n}[n;{\bf R}]^{a_1\cdots a_n}_{1\cdots
n}\omega^2\cdots}
where the descendents in $b_{-2}(\omega^1)$ anihilate against $\omega^0$ to the
left, and so do not contribute in the
first line. We move the
$\theta^{-1}$ mode to the left in the second line, picking up powers of $P{\bf
R}$. The third equality then applies $b_2$. Only
its action on $\theta^{-1}$ contributes, since modes $\theta^2$ or higher can
be moved to the right and anihilate. The fourth equality
moves the resulting $\theta^1$ to the right, picking up powers of $P{\bf R}$
again.

We now use the Hecke condition in the form
$(P{\bf R})^2=q^{-2}+(q^{-2}-1)P{\bf R}$ and the Yang-Baxter equations in the
form $(P{\bf R})_{23}(P{\bf R})_{12}(P{\bf R})_{23}=
(P{\bf R})_{12}(P{\bf R})_{23}(P{\bf R})_{12}$ repeatedly, to observe that
\align{&&\nquad \theta^0\cdots\theta^0\theta^1\overline{[1,n+1;{\bf
R}]}[1,n+1;{\bf R}]\\
&&=\theta^0\cdots\theta^0\theta^1( q^{-2}\overline{[2,n+1;{\bf R}]}[2,n+1;{\bf
R}]\\
&&\quad +(q^{-2}-1)(P{\bf R})_{nn+1}\cdots(P{\bf R})_{23}(P{\bf R})_{12}(P{\bf
R})_{23}\cdots (P{\bf R})_{nn+1})\\
&&=\theta^0\cdots\theta^0\theta^1\left( q^{-2}\overline{[2,n+1;{\bf
R}]}[2,n+1;{\bf R}] +(q^{-2}-1)[1,n;{\bf R}]
(P{\bf R})_{nn+1} \overline{[1,n;{\bf R}]}\right)\\
&&=\theta^0\cdots\theta^0\theta^1\left( q^{-2}\overline{[2,n+1;{\bf
R}]}[2,n+1;{\bf R}]
+(q^{-2}-1) q^{-2(n-1)} \overline{[1,n+1;{\bf R}]}\right)\\
&&=\cdots=\theta^0\cdots\theta^0\theta^1
\left(q^{-2n}+(q^{-2}-1)q^{-2(n-1)}(\overline{[n+1;{\bf R}]}-1\right)\\
&&=\theta^0\cdots\theta^0\theta^1\left(q^{-2n}-(q^{-2}-1)q^{-2(n-1)}\right)
=\theta^0\cdots\theta^0\theta^1q^{-2(n-1)}.}
The third equality replaces $P{\bf R}$ by $q^{-2}$ in $[1,n;{\bf R}]$ since it
acts on $\theta^0\cdots\theta^0$ to its left.
We then iterate these steps, collecting the $\overline{[\ ,n+1 ;{\bf R}]}$
which are  generated in this way as
$\overline{[n+1;{\bf R}]}-1$. Finally, we note that
\[ \theta^0\cdots\theta^0\theta^1\overline{[n+1;{\bf
R}]}=\theta^0\cdots\theta^0\overleftarrow{\del}\cdot\theta^1=0\]
since on the right hand side we have the braided differential of $n+1$ copies
of $\theta^0$, which vanishes.

With this result, we can complete our calculation as
\[b_2(\omega^0b_{-2}(\omega^1))\omega^2\cdots=\omega^0q^{-2(n-1)}
\theta_{a_1}^1\cdots\theta^1_{a_n}[n;{\bf R}]^{a_1\cdots a_n}_{1\cdots
n}\omega^2\cdots
=q^{-2(n-1)}[n;q^{-2}]\omega\]
since $P{\bf R}$ can be replaced by $q^{-2}$ when acting on the algebra
$\Lambda^{(1)}(R)$.

Adding this contribution to that from Lemma~3.3, we find
\[ b_2(b_{-2}(\omega))=\left([n;q^{-2}](1+q^{-2(n-1)})
+(1-q^{-2})([n-1;q^{-4}]-q^{-2(n-1)}[n-1;q^{-2}])\right)\omega\]
which computes to the final result stated. \endproof

Although we have only covered the $i=1,2$ cases of (\ref{heis}) in this paper,
it is clear that the method
introduced here can provide a viable alternative to the vertex operator proof
in \cite{KasMiwSte}. Since
the approach there uses directly the correlation function for $XXZ$ vertex
operators, our direct `braided geometric'
technique implies in principle a new approach the the computation of these.

\section{Concluding remarks}

It is significant that all computations in this paper have been made without
reference to any specific details of the $R$-matrix,
so long as it is Hecke type. This means that the fermionic Fock space
construction in \cite{KasMiwSte} works quite
generally; it may be interesting to consider some non-standard examples. A
further question is
how to extend the above methods to non-Hecke cases such as the affine quantum
group $U_q(\hat{so_3})$. Related to this is
the construction of higher level fermionic Fock space representations, even for
$U_q(\hat{sl_2})$. For these one should
make semi-infinite tensor products of fermionic quantum planes where the
underlying finite-dimensional
R-matrix is not of Hecke type. We note that the Baxterisation formula for the
parametrised $R$-matrix in the $\hat{so_n}$ case
is indeed known, though having now a more complicated form. Hence in principle
our `decomposition' methods might be applied.

Also, in braided geometry the fermionic quantum planes (like
other quantum planes) are fully covariant not exactly  under $U_q(sl_n)$ (or
other quantum group, according to the R-matrix) but
under a dilatonic extension of it. This is needed whenever the quantum plane
normalisation is not the quantum group normalisation
of the R-matrix. Analogously, the fermionic Fock space is not quite covariant
under the quantum loop group associated to
$R(z)$ but under its central extension, which in our case is $U_q(\hat{sl_n})$.
Formally, and before considering the
R-matrix normalisation, the exchange algebra (\ref{fock}) should be covariant
under the quantum loop group in the R-matrix
form with generators $\vecl^\pm(z)$, which would make it a level 0 module of
$U_q(\hat{sl_n})$. Hence it appears that
similar `dilaton' effects are responsible for the
anomaly which makes the fermionic Fock space
considered above into level 1. This is another direction for further work.

\end{document}